\documentclass[draftcls,12pt, onecolumn]{IEEEtran}
\usepackage{amsmath,amssymb}
\usepackage{ifpdf}
\usepackage{subfigure}
\usepackage{indentfirst}
\usepackage{url}
\usepackage{times}
\usepackage{cite}
\usepackage[body={17true cm,22true cm},top=2.5cm]{geometry}

\ifpdf
  \usepackage[pdftex]{graphicx}
  \usepackage[pdftex,CJKbookmarks=true,bookmarksnumbered=true,bookmarksopen=true,
     colorlinks=true,citecolor=magenta,linkcolor=red,anchorcolor=green,urlcolor=blue]{hyperref}
   \DeclareGraphicsExtensions{.pdf}
\else
   \usepackage[dvips]{graphicx}
   \usepackage[dvipdfm,CJKbookmarks=true,bookmarksnumbered=true,bookmarksopen=true,
     colorlinks=true,citecolor=blue,linkcolor=red,anchorcolor=green,urlcolor=blue]{hyperref}
   \DeclareGraphicsExtensions{.eps}
\fi

\newcommand{\complex}{\mathcal{C}}
\newcommand{\real}{\mathcal{R}}

\graphicspath{{../Pic/}} \DeclareGraphicsExtensions{.pdf,.epsf,.png}

\DeclareMathOperator{\rank}{rank} \DeclareMathOperator{\diag}{diag}

\newtheorem{theorem}{Theorem}
\newtheorem{lemma}{Lemma}
\newtheorem{corollary}{Corollary}
\newtheorem{remark}{Remark}

\allowdisplaybreaks

\numberwithin{equation}{section}

\title{Stability and steady state analysis of distributed cooperative droop controlled DC microgrids
\thanks{This research was supported by a grant from the National 863 Program of China (2015AA050104) and the National Natural Science Foundation of China (61074122, 61104149).}
}

\author{Ji~Xiang, 
        Yu~Wang, Yanjun~Li, and
        Wei~Wei,
\thanks{J. Xiang, Y. Wang and W. Wei are with the Department of System Science and
Engineering, College of Electrical Engineering, Zhejiang University,
P. R. China. (e-mail: jxiang@zju.edu.cn)}
\thanks{Y. Li is with the School of Information and Electrical Engineering, Zhejiang University City College, P. R. China. (e-mail: liyanjun@zucc.edu.cn.)}
}

\begin{document}
\maketitle %
\begin{abstract}
Distributed cooperative droop control consisting of the primary decentralized droop control and the  {secondary} distributed correction control  {is studied in this paper, which aims to achieve an exact current sharing between generators, worked in the voltage control mode, of DC microgrids.} For the DC microgrids with the distributed cooperative droop control, the dynamic stability has not been well investigated although its steady performance has been widely reported. This paper focuses on the stability problem and shows it is equivalent to the semistability problem of a class of second-order matrix systems. Some further sufficient conditions as well followed. The steady state is analyzed deeply for some special cases. A DC microgrid of three nodes is simulated on the Matlab/Simulink platform to illustrate the efficacy of analytic results.
\end{abstract}

\section{Introduction}
As growing stress on the environment protection and the depletion of fossil energy, renewable energy sources play the role more and more important in the power generation, such as wind and photovoltaic (PV) \cite{DincerRSER2000}. These sources output power by distributed generators (DGs) that are usually integrated to the distribution system by  {microgrids}. A microgrid is a collection of DGs, distributed storage systems and loads, which is connected to the grid via a point of common coupling (PCC) \cite{HatziargyriouPEM2007, YazdanianTSG2014}. 

Since most  {renewable energy sources (RES)} and storage systems, as well as many loads like vehicles, data centers, and telecom systems, have a nature of direct current, it is  {more preferable to connect} these RESs, storage systems and the load to form a DC microgrid by using dc-dc converters directly without two stages dc-ac-dc conversion. Only one dc-ac converter is applied at PCC  {for a connection} to the grid. There are no reactive power and frequency synchronization in a DC microgrid, both of them being the main challenging problems of AC systems. Owing to these striking features, DC microgrids have been attracting considerable attention in more recent years \cite{JustoRSER2013, ElsayedEPSR2015, PlanasRSER2015}.

Similar with AC microgrids, an important control objective of DC microgrids is to share the power demanded by loads among different sources. The sharing control for microgrids including both AC and DC can be roughly categorized into three types: centralized, decentralized and distributed  {control} \cite{MokhtariTPS2013}.  {The centralized control} coordinates all sources in an optimal way but requires a high-bandwidth communication to timely collect all the information so as to have low reliability and expandability.  {The decentralized control} often adopts a droop control to make a sharing among sources. But the sharing may be not effective due to a lack of broader available information. For a DC microgrid particularly there is not a parameter which remains the same throughout the DC microgrid, like the frequency for AC microgrids \cite{AnandTPE2013}.  {The distributed control} as a strategy between the centralized  {control} and the decentralized  {control} is more robust and expandable. It can make an exact sharing among all sources in the cost of low-bandwidth communications between neighboring sources.   

DC decentralized droop control can be traced back to about twenty years ago for a superconductive DC system \cite{JohnsonTPD1993}, where the voltage is uniform for all terminals and thereby is utilized as the same parameter like  {the} frequency of AC system{s}, to coordinate terminal currents. In the presence of line resistances, the DC voltage is no longer  {the} uniform measure such that the load sharing is difficult to obtain by the decentralized droop control.  In \cite{HaileselassieTPS2012}, it has been illustrated that with line resistances the power sharing has a large deviation from that in a lossless DC microgrid. Such a sharing deviation can be reduced by large droop gains which however  {lead} to a large deviation of terminal voltage. In \cite{AnandTPE2013}, the two factors hampering the application of small droop gains are presented. In \cite{BeertenTPS2013}, the steady state of decentralized droop control is addressed. The influences of droop gains on the power sharing and voltage deviations are illustrated and an optimal droop gain setting problem is issued.

To improve the current sharing accuracy, a hierarchical controller is presented in \cite{GuerreroTIE2011} where the large droop gains are adopted in the primary control for a small deviation current sharing and the large deviation of voltage is compensated by the secondary control. In \cite{AnandTPE2013}, a distributed droop controller based on the average current of all terminals is proposed to make an exact current sharing. But an extra wire laying along with power lines is required to connect the measured values of currents such that the information of average current is available in real time. In \cite{LuTPE2014}, a distributed droop control based on a low-bandwidth communication is proposed, in which the converter's voltage and current are exchanged between neighboring converters. The stability analysis has been made only on a two-node DC microgrid. A distributed droop controller including two modules, voltage regulator and current regulator for a meshed DC microgrid has been proposed in \cite{NasirianTPE2015}. An extension to adaptive droop gains is reported in \cite{NasirianTEC2014}. Although the analysis of steady states and performances has been made in the frequency domain, the stability of overall closed-loop system is not addressed.    

It is difficult to analyze the overall stability of microgrids. The conventional frequency domain method is difficult to tackle the whole dynamics with high dimensions.  A few works on this topic often are based on the small-signal model. In \cite{DragicevicTPE2014}, the stability of battery converter with adaptive droop gains has been addressed. In \cite{RadwanTSG2012}, it was shown that even if each converter is stable by itself the stability of overall DC microgrid is not ensured because of the coupling between converter regulators. In \cite{AnandTIE2013}, a linearized model including sources, lines and loads is presented for DC microgrids. The eigenvalues of the system matrix determine the stability of DC microgrid. The relationship between the eigenvalue locations and the line impedance  {is} discussed. 

This paper addresses the distributed cooperative droop controller with only one module of current sharing regulator. The controller has two control levels. The primary control is  {a} decentralized voltage droop control to regulate the converter output voltage according to its output current. The secondary control is  {a} distributed current sharing control in the sense that the neighboring converters exchange their p.u. currents via a low-bandwidth communication. Our focus is on the stability problem owing to the secondary distributed control. The primary droop control is often performed at the DC/DC or AC/DC converters that have a fast response compared to the secondary distributed control. Its dynamics therefore is omitted, as well as the effect of linear inductances that are very small in DC microgrids. Our goal is to find the influences of droop gains, controller gains and admittance matrix on the stability of overall system. 

 {As for as the current sharing is concerned, the control goal is to force the currents of generator to reach the same value (or the same ratio with respect to their maximal/nominal currents). The desired sharing currents can not be assigned a prior because of the fluctuation of loads in a power system. Such a scenario can be characterized by the notation of semistablity \cite{CampbellSIAMJMA1979, BernsteinJMD1995, BhatACC1999}, which means that the steady state is not completely determined by the system dynamics, but depends on the system initial conditions as well. Semistability is an appropriate notation for the analysis of self-organized behaviors of networked systems which rely on the initial configuration, and has been applied for the consensus problem of linear \cite{HuiIJC2009, HuiSCL2011} and nonlinear networked systems \cite{HuiAUT2008, HuiTAC2008}. This paper will also use the tool of semistability to analysis the stability of the closed-loop system under the secondary distributed current sharing control. 
More recently, Andreasson et. al., regarded the terminals of HVDC transmission systems as the controlled current sources and presented three kinds of distributed controllers to regulate the terminal voltages, as well as the related sufficient conditions for the stability of the closed-loop system in \cite{AndreassonARXIV2014}. Zhao et. al., also investigated a two-level current sharing control problem of a class of DC microgrids consisting of current sources and constant current loads in \cite{ZhaoACC2015sub}. A decentralized droop control that can achieve the current sharing or optimal economic dispatch with suitable gains is presented. Then a secondary distributed control is presented to compensate the voltage drifts due to the primary decentralized droop control. } Both of the two recent references consider the model in which the generator current is to be controllable. Considering that many practical DC/DC converters work as voltage sources, this paper addresses the current sharing problem for the case of the voltage to be controlled. Moreover, the stability analysis here goes along the line of the semistability which is as well different from the line of the characteristic equation used in \cite{AndreassonARXIV2014}.

The remainder of this paper is structured as follows. The dynamic model of the decentralized droop control is presented in Section \ref{sec02}, where the steady state and stability  under the primary decentralized droop control are detailed. The distributed current sharing control is presented in Section \ref{sec03}. The semistability of the overall system is addressed since the current sharing is not a decaying requirement. The steady state is detailed as well. A simulation example is illustrated in Section \ref{sec04} to validate the analytic results, followed by a conclusion in Section \ref{sec05}.

\section{Decentralized droop control} \label{sec02}
\subsection{Dynamic model}
Consider a DC microgrid of $n$ generator nodes and $n_l$ constant impedance load nodes. Without loss of generality assume the first $n$ nodes are generators. Denote two node sets by $\mathcal{N}=\{1,2,\cdots,n\}$ and $\mathcal{N}_l=\{n+1,\cdots,n+n_l\}$.  A line connecting node $k$ and $j$ is associated with a branch conductance $G'_{kj}\ge 0$. $G'_{kj}=0$ if and only if there is no connection between node $k$ and $j$. A node $k$ is associated with a current injection $i_k$, an output voltage $u_k$ and a shunt conductance $G'_{kk}\ge 0$. For all $k\in\mathcal{N}_l$, $G'_{kk}>0$. If $G'_{kk}=0$, then node $k$ has no local load. There is no current injection for load nodes, so $i_k=0$ for all $k\in\mathcal{N}_l$.

\begin{figure}
\centering
\includegraphics[width=0.4\textwidth]{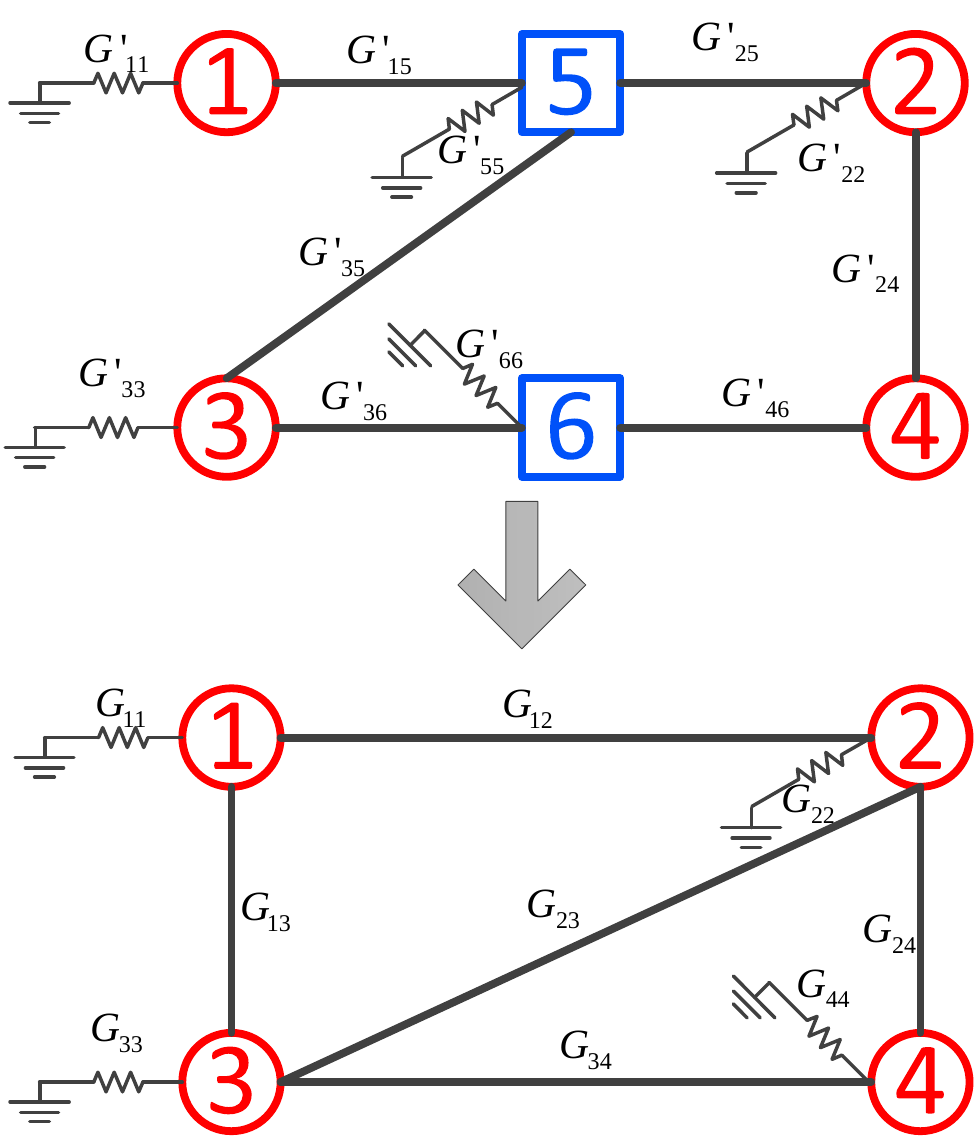}\\
\caption{A DC microgrid with $n=4$ generator nodes, $n_l=2$ load nodes and $m=6$ lines. Generator $4$ has no shunt conductance. After a network reduction, it becomes one with only $4$ generator nodes and $5$ lines. Each node has a shunt conductance connected.}\label{figa01}
\end{figure}

Define a stacked current vector $I'\in\real^{n\times 1}$ by $I'=\mathrm{col}(i_1,i_2,\cdots,i_n)$ and stacked voltage vectors $U'\in\real^{n\times 1}$ and $U'_L\in\real^{n_l\times 1}$ by $U'=\mathrm{col}(u_1,\cdots,u_n)$ and $U'_L=\mathrm{col}(u_{n+1},\cdots,u_{n+n_l})$, respectively. Then the DC microgrid current-voltage equation is given by
\begin{equation} \label{eqII01}
\begin{bmatrix} I' \\ 0 \end{bmatrix} = Y'\begin{bmatrix} U' \\ U'_L \end{bmatrix}=
\begin{bmatrix} Y'_{gg} & Y'_{gl} \\ Y'_{lg}  & Y'_{ll} \end{bmatrix}
\begin{bmatrix} U' \\ U'_L \end{bmatrix} 
\end{equation}
where the conductance matrix $Y'=[Y'_{kj}]\in\real^{(n+n_l)\times (n+n_l)}$ is defined as
\begin{equation}
Y'_{kj}=\left\{\begin{array}{ll} \sum_{j=1}^nG'_{kj} & k=j \\ -G'_{kj} & k\neq j
\end{array},\right.
\end{equation}
$Y'_{gg}\in\real^{n\times n}$, $Y'_{ll}\in\real^{n_l\times n_l}$ and $Y'_{gl}=(Y'_{lg})^T$. Recalling $G_{kk}>0$ for all $k\in\mathcal{N}_l$, $Y'_{ll}$ is a symmetric positive definite (SPD) matrix. Solving $U'_L$ from the last $n_l$ equation of \eqref{eqII01} yields
\begin{equation}
I'=(Y'_{gg}-Y'_{gl}Y_{ll}^{-1}(Y'_{gl})^T)U',
\end{equation}
by which a network reduction can be carried out to obtain a new DC microgrid with only  generator nodes. An example is given in Fig. \ref{figa01}. 

Hereafter we, without loss of generality, consider a DC microgrid with $n$ generator nodes, which is associated with current vector $I$, voltage vector $U$ and conductance matrix $Y=Y'_{gg}-Y'_{gl}Y_{ll}^{-1}(Y'_{gl})^T=[Y_{kj}]\in\real^{n\times n}$,  {satisfying
\begin{equation} \label{eqa01}
I=YU.
\end{equation}}

 {
\begin{remark}
The network reduction is often used in the analysis of electrical networks  and power systems \cite{DorflerTCASI2013}, which simplifies the loads as the equivalent shunt loads and the branch impedances such that the influences of generators are focused. The developed above is also feasible for the constant current loads, for which the current injection of the load nodes is not zero but $I'_L$. A similar operation yields $I'-Y'_{gl}(Y'_{ll})^{-1}I'_L=(Y'_{gg}-Y'_{gl}Y_{ll}^{-1}(Y'_{gl})^T)U'$. Taking $I=I'-Y'_{gl}(Y'_{ll})^{-1}I'_L$ obtains \eqref{eqa01} as well. 
\end{remark}}

\bigskip
Generators in the microgrid connects to the DC bus via DC/DC converters that in general have a very small time constant. Due to the ultra fast responses of converters, the generator can be simplified as a DC voltage source whose voltage is regulated instantaneously, namely,
\begin{equation} \label{eqII04}
u_k=u_k^{ref},\quad k\in\mathcal{N},
\end{equation}
where $u_k$ and $u^{ref}_k$ are the output voltage and the output voltage reference of node $k$. A DC microgrid makes a load current sharing by a  {voltage droop} controller with which the voltage reference of generator will reduce when its output current increases, being implemented as 
\begin{equation} \label{eqII05}
u_k^{ref}=u_k^d-R_ki^m_k, 
\end{equation}
where $u_k^d$ is the rated output voltage, $R_k$ is the internal resistance that might be a virtual one to be designed and $i^m_k$ the output current measured. 

Since of interest is the DC current, a low pass filter on the output current $i_k$ is used to get $i_k^m$,
\begin{equation} \label{eqII06}
\tau_k \dot{i}_k^m = -i_k^m+i_k,
\end{equation}
where $\tau_k$ is the time constant of low-pass filter, with which the cut-off frequency is given by $1/\tau_k$.  {According to \eqref{eqa01}, one has}
\begin{equation} \label{eqII07}
i_k= \sum_{j=1}^n Y_{kj}u_j, \quad k\in\mathcal{N}.
\end{equation}
Combining \eqref{eqII04} and \eqref{eqII05} and replacing \eqref{eqII07} into \eqref{eqII06} yield the following dynamic equation of  {the} DC microgrid,
\begin{equation} \label{eqII08}
\tau_k\dot{i}_k^m=-i_k^m+\sum_j Y_{kj}(u_j^d-R_ji_j^m), \quad k\in\mathcal{N}.
\end{equation}
Define $I^m=\mathrm{col}(i_1^m,\cdots,i_n^m)$ and $U^d=\mathrm{col}(u_1^d,\cdots,u_n^d)$, then the compact form of \eqref{eqII08} is given by
\begin{equation} \label{eqII09}
D\dot{I}^m = -(E+YR)I^m+YU^d,
\end{equation}
where $D=\mathrm{diag}(\tau_1,\tau_2,\cdots,\tau_n)$, $R=\mathrm{diag}(R_1,R_2,\cdots R_n)$ and $E$ denotes the identity matrix. Equation \eqref{eqII09} describes the dynamics of  {the} DC microgrid under droop controller \eqref{eqII05} that depends on local measurements and therefore is decentralized.

\subsection{Stability analysis}
Before presenting the stability result of \eqref{eqII09}, the shunt conductance of each node is separated from the conductance matrix $Y$ by $Y=Y_s+Y_c$, where $Y_s=\mathrm{diag}(G_{11},G_{22},\cdots,G_{nn})$ denotes the shunt conductance and $Y_c$ consists of branch conductances, defined by
\begin{equation}
Y_c=\begin{bmatrix}
\sum\limits_{j\neq 1} G_{1j} & -G_{12} & \cdots & -G_{1n} \\
-G_{21} &   \sum\limits_{j\neq 2} G_{2j} & \cdots & -G_{2n}\\
\vdots & \vdots & \ddots & \vdots\\
-G_{n1} & -G_{n2} & \cdots & \sum\limits_{j\neq n} G_{nj}
\end{bmatrix}.
\end{equation}
\begin{lemma}
Given $\tau_k>0$ and $R_k>0$ for all $k\in\mathcal{N}$, the DC microgrid given by \eqref{eqII06} and \eqref{eqII07} is exponentially stable with decay rate not larger than $-\psi$ under droop controller \eqref{eqII05},  where $\psi=\min\limits_{i\in\mathcal{N}} (\tau^{-1}_i+\tau_{i}^{-1}R_iG_{ii})$ 
\end{lemma}
\begin{IEEEproof}
With controller \eqref{eqII05}, the DC microgrid has a closed-loop dynamic described by \eqref{eqII09} whose system matrix is $-D^{-1}(E+YR)$. Taking a similar transformation obtains 
\begin{equation}
R(D^{-1}(E+YR))R^{-1}=D^{-1}+D^{-1}RY
=D^{-1}+D^{-1}RY_s+D^{-1}RY_c\triangleq \Psi_1.
\end{equation}
Since $Y_c$ has all rows of zero sum, the left bound of all eigenvalues of $\Psi_1$ is not less than $\frac{1+G_{ii}R_{i}}{\tau_{i}}$ for all $i\in\mathcal{N}$ by the Gershgorin circle theorem. Thus, the proof is completed.
\end{IEEEproof}

\begin{remark}
In general, the shunt conductance $G_{ii}$ is small, as well as the virtual resistance $R_i$, therefore it is the maximal time constant $\tau_{i}$ that determines the decay rate.
\end{remark}
  
\subsection{Steady state}
Since $\tau_i>0$ for all $i\in\mathcal{N}$, the steady state can be solved from setting the right side of \eqref{eqII09} to zero,
\begin{equation} \label{eqII12}
I^{ss}=I^{m^{ss}}=(E+YR)^{-1}YU^d,
\end{equation}
where the superscript $ss$ denotes the steady state of variable and the first equality comes from \eqref{eqII06}. Subsequently, the steady output voltage of each node is given by
\begin{equation} \label{eqII13}
U^{ss}=U^d-RI^{ss}=(E+RY)^{-1}U^d.
\end{equation}
Denote by $\underline{1}_n$ a n-order vector with all elements being $1$. Define the operator measuring the deviation from the current sharing by $\sigma = \|(E-\frac{1}{n}\underline{1}_n\underline{1}_n^T)I^{ss}\|$. Then the following result can be obtained,

\begin{lemma} \label{th02}
The steady state \eqref{eqII12} has the following properties:
\begin{enumerate}
\item Given $U^d=u^d\underline{1}_n$, $R=rE$ and $Y_s=gE$, then all the nodes have the same steady state current, $I^{ss}=\frac{g}{1+gr}u_d\underline{1}_n$.
\item Given $U^d=u^d\underline{1}_n$, then $\dfrac{\sigma^2}{u_d^2}\le n\left(\frac{\bar{g}}{1+\bar{g}\underline{r}}\right)^2-\frac{1}{n}\left(\sum\limits_{i=1}^n \frac{G_{ii}}{1+G_{ii}R_i}\right)^2$, where $\bar{g}=\max_i G_{ii}$ and $\underline{r}=\min_i R_i$.
\end{enumerate}
\end{lemma}
\begin{IEEEproof}
The first property is clear and its proof is omitted. Below we show the second property. 
Noting $I^{ss}=(Y^{-1}+R)^{-1}U^d$ and $(E-\frac{1}{n}\underline{1}_n\underline{1}_n^T)^2=(E-\frac{1}{n}\underline{1}_n\underline{1}_n^T)$, it follows that
\begin{equation}
\frac{\sigma^2}{u_d^2}=\underline{1}_n^T(Y^{-1}+R)^{-2}\underline{1}_n-\frac{1}{n}(\underline{1}_n^T(Y^{-1}+R)^{-1}\underline{1}_n)^2.
\end{equation}
For a connected network, there is a transformation matrix $T=[\underline{1}_n,\star]$  diagonalizing $Y_c$ such that $T^{-1}Y_cT=\mathrm{diag}(0,\lambda_2, \cdots,\lambda_n)$,
where $T^{-1}$ has the form of $T^{-1}=[\frac{1}{n}\underline{1}_n,\star]^T$ and $\star$ denotes the part of no interest. It can be verified that $T^{-1}\underline{1}_n=[1,0,\cdots,0]^{T}$ and $\underline{1}_n^T T=[n,0,\cdots,0]$. Notice 
\begin{equation}
Y^{-1}+R\ge (\bar{g}E+Y_c)^{-1}+\underline{r}E= T^{-1}\mathrm{diag}(\bar{g}^{-1}+\underline{r},
(\bar{g}+\lambda_2)^{-1} +\underline{r},\cdots,(\bar{g}+\lambda_n)^{-1}+\underline{r})T,
\end{equation}
and 
\begin{equation}
Y^{-1}+R\le Y_s^{-1}+R.
\end{equation}
The following inequality holds,
\begin{multline}
\frac{\sigma^2}{u_d^2}\le \underline{1}_n^TT \big[\mathrm{diag}(\bar{g}^{-1}+\underline{r},(\bar{g}+\lambda_2)^{-1}+\underline{r},\cdots,
(\bar{g}+\lambda_n)^{-1}+\underline{r})\big]^2T^{-1}\underline{1}_n\\-\frac{1}{n}\big[\underline{1}_n^T\mathrm{diag}\big((G_{11}^{-1}+R_1)^{-1},\cdots,(G_{nn}^{-1}+R_n)^{-1}\big)\underline{1}_n \big]^2,
\end{multline}
which completes the proof.
\end{IEEEproof}

\medskip
With the condition of the first property, $G_{ii}=\bar{g}$ and $R_i=\underline{r}$ for all $i\in\mathcal{N}$,  $n\left(\frac{\bar{g}}{1+\bar{g}\underline{r}}\right)^2-\frac{1}{n}\big(\sum_{i=1}^n \frac{G_{ii}}{1+G_{ii}R_i}\big)^2=0$ and therefore $\sigma=0$, achieving the same current between nodes. 

 {In the case that} $G_{ii}$ and $r_i$ is much less than $1$, so approximately $\frac{\sigma^2}{u_d^2}\le n\big((\bar{g})^2-(\tilde{g})^2\big)$ where $\tilde{g}=\frac{1}{n}\sum_i G_{ii}$ denotes the average shunt conductance. For a large virtual resistance satisfying $R_iG_{ii}\gg 1$, approximately $\frac{\sigma^2}{u_d^2}\le n\Big((\underline{r}^{-1})^{2}-(\frac{\sum_i R_i^{-1}}{n})^2\Big)$.

\bigskip
Lemma 2 discusses two cases about the current sharing under the decentralized droop control when generators have the same nominal output voltage. The first property shows that if all the generators have the same shunt impedance and virtual resistance then  they have the same currents regardless of the branch impedances. The second property shows that the nominal deviation of current sharing converges to zero as the virtual resistance $R_i$ converge to infinity and is determined by the difference of the shunt impedances if $R_i$ is sufficiently small.

\section{Distributed current sharing control} \label{sec03}
As stated in the above section, the same voltage reference $U^d=u^d\underline{1}_n$ in general  can not lead to the same current. This section address a distributed method to realize a current sharing. Only the neighboring nodes exchange information each other and each node adjusts the rated voltage according to the current bias from their neighboring nodes. A low-bandwidth communication channel is required between the neighboring nodes.

Let $\mathcal{G}=\{\mathcal{N},\mathcal{E}\}$ express the information flow between nodes, where $\mathcal{E}\subseteq \mathcal{N}\times \mathcal{N}$ denotes the edge set. Let $L=(l_{ij})\in\real^{n\times n}$ be the associated Laplacian matrix of $\mathcal{G}$. If $Y_{ij}\neq 0$, then $l_{ij}=-1$; or else $l_{ij}=0$. The information graph has no self-loop, so $l_{ii}=\sum_{j=1,j\neq i}^n l_{ij}$.

In contrast to the same current, the same current ratio is a more rational index for current sharing. Denote by $I_i^{max}$ the maximum current of node $i$, then our goal is to realize
\begin{equation}
\frac{i_1}{I_1^{max}}=\frac{i_2}{I_2^{max}}=\cdots=\frac{i_n}{I_n^{max}}.
\end{equation}
Let current ratio $i_k^r=\frac{i^m_k}{I_k^{max}}$ be the information exchanged between neighboring nodes. The following PI controller is proposed for the rated output voltage
\begin{equation} \label{III02}
u_k^{d}=(-\alpha_k-\frac{\beta_k}{s})\sum_j l_{kj}i^r_{j}, \quad k\in\mathcal{N},
\end{equation}
where $\alpha_k>0$ and $\beta_k>0$ are proportional and integrator gains of node $k$, respectively. 
Replacing \eqref{III02} into \eqref{eqII08}, the closed-loop dynamics of each node is described by
\begin{equation} \label{III03}
\left\{
\begin{split}
&\dot{i}^m_k=-\frac{1}{\tau_k} i_k^m+\frac{1}{\tau_k}\sum_{j=1}^n Y_{kj}u_j^d-\frac{1}{\tau_k}\sum_{j=1}^{n}Y_{kj}R_ji_j^m\\
&\dot{u}_k^d=-\alpha_k\sum_{j=1}^n l_{kj}\frac{\dot{i}_j^m}{I_j^{max}}-\beta_k \sum_{j=1}^n l_{kj}\frac{i_j^m}{I_j^{max}}
\end{split}
\right.,\,\,\,k\in\mathcal{N}.
\end{equation} 
Define $\Upsilon=\mathrm{diag}(I_1^{max}, \cdots,I_n^{max})$, $\Phi=\mathrm{diag}(\alpha_1,\cdots,\alpha_n)$ and $\Psi=\mathrm{diag}(\beta_1,\cdots,\beta_n)$. Then the above closed-loop system can be rewritten as the following compact form,
\begin{equation} \label{III04}
\begin{bmatrix} \dot{I}^m \\ \dot{U}^d \end{bmatrix}=
\begin{bmatrix}
-D^{-1}(E+YR) & D^{-1}Y \\
M_{21} & M_{22}
\end{bmatrix}
\begin{bmatrix} I^m \\ U^d \end{bmatrix},
\end{equation}
where $M_{21}=-\Psi L (\Upsilon)^{-1}+\Phi L (\Upsilon)^{-1} D^{-1}(E+YR)$ and $M_{22}=-\Phi L (\Upsilon)^{-1}D^{-1}Y$.

 {
\begin{remark}
We would like to point out that the secondary control \eqref{III02} aims to make the current sharing between nodes and is distributed because of the use of the information of neighboring nodes. While the secondary control in \cite{GuerreroTIE2011} aims to compensate the voltage drifts due to the primary droop control and is decentralized because only the local voltage and reference voltage are used.  
\end{remark}}

\subsection{Stability analysis}
As a linear autonomous system, system \eqref{III04} is required to be neither unstable nor asymptotically stable. An asymptotical stability means output voltages and currents converging to zero,  {which certainly is not what we want}. Actually, we hope that the nodes have the same output current ratio with their output voltage in the desired region. Such a case corresponds to the term of {\em semistable} that is recalled in the Appendix. 

\begin{lemma} \label{th03}
System \eqref{III04} is semistable if and only if the following matrix second-order system 
\begin{equation} \label{III05}
D\Upsilon\ddot{x}(t)+(\Upsilon+YR\Upsilon+Y\Phi L)\dot{x} + Y\Psi L =0,
\end{equation}
is semistable.
\end{lemma}
\begin{IEEEproof}
Let $A_c$ denote the system matrix of \eqref{III04}, which is equivalent to 
\begin{equation}
A_c=
\begin{bmatrix} E & 0 \\ -\Phi L (\Upsilon)^{-1} & E \end{bmatrix}
\begin{bmatrix} -D^{-1}(E+YR) & D^{-1}Y \\ -\Psi L (\Upsilon)^{-1} & 0 \end{bmatrix}.
\end{equation}
Notice that matrix $A_c$ is similar to 
\begin{equation} \label{III07}
\tilde{A}_c = \begin{bmatrix} -D^{-1}(E+YR+Y\Phi L(\Upsilon)^{-1}) & D^{-1}Y \\ -\Psi L (\Upsilon)^{-1} & 0 \end{bmatrix},
\end{equation}
which satisfies
\begin{equation}
\begin{bmatrix}
\lambda_i+D^{-1}(E+YR+Y\Phi L(\Upsilon)^{-1}) & -D^{-1}Y \\ \Psi L (\Upsilon)^{-1} & \lambda_i
\end{bmatrix}\begin{bmatrix} x_{i1} \\ x_{i2} \end{bmatrix}
=0
\end{equation}
where $\lambda_i$ is the $i^{th}$ eigenvalue of $\tilde{A}_c$ and $x_i=[x_{i1}^T,x_{i2}^T]^T\in\complex^{2n}$ is the corresponding complex eigenvector. Noting that $-\Psi L (\Upsilon)^{-1}x_{i1}=\lambda_i x_{i2}$ and subsequently $-\lambda_i D^{-1}Yx_{i2}=D^{-1}Y\Psi L (\Upsilon)^{-1}x_{i1}$, it follows that
\begin{equation}
\big(\lambda_i^2+ D^{-1}(E+YR+Y\Phi L(\Upsilon)^{-1})\lambda_i 
+D^{-1}Y\Psi L(\Upsilon)^{-1}\big)x_{i1}=0
\end{equation}
Multiplying $D$ on both sides and defining $\tilde{x}_{i1}=(\Upsilon)^{-1}x_{i1}$ yield
\begin{equation}
\big(D\Upsilon\lambda_i^2+ (\Upsilon+YR\Upsilon+Y\Phi L)\lambda_i+Y\Psi L\big)\tilde{x}_{i1}=0
\end{equation}
The term in bracket is just the characteristic equation of dynamic system \eqref{III05}. This completes the proof.
\end{IEEEproof}

\medskip
 {
The equation \eqref{III07} shows that the system \eqref{III04} can be rewritten as 
\begin{eqnarray*}
D\dot{I}^m&=&-(E+YR+Y\Phi L(\Upsilon)^{-1})I^m+Y\Theta, \\
\dot{\Theta}&=&-\Psi L \Upsilon^{-1} I^m,
\end{eqnarray*}
with $\Theta=U^d-\Phi L\Upsilon^{-1}I^m$. It is seen that $\Theta$ changes and then forces $I^m$ to change until the current sharing has been achieved. That is, $L\Upsilon^{-1} I^m=0$ and subsequently $\Theta=U^d$.}

Lemma \ref{th03} converts the stability of \eqref{III04} into a stability problem of a matrix second-order system that by itself owns fundamental importances in many fields, such as vibration and structure analysis, spacecraft control and robotics control. There are many results for the stability analysis of \eqref{III05} with symmetric matrix coefficients, which however are not directly applied here due to asymmetric matrix coefficients arising from the heterogeneity between nodes. Below we further present some sufficient condition for the stability of \eqref{III04} for some special cases. Before to proceed, the following results are recalled \cite{BernsteinJMD1995}. 

Given a second-order dynamic system $(M_s+M_{k})\ddot{x}(t)+(D_s+D_k)\dot{x}+(K_s+K_k)=0$, its eigensolution can be written as
\begin{equation} \label{III11}
\big(\lambda_i^2(M_s+M_k)+\lambda_i(D_s+D_k)+(K_s+K_k)\big)x_i 
 =0,\quad i=1,\cdots,2n,
\end{equation}
where $\lambda_i$ and $x_i$ are the $i$th eigenvalue and the corresponding complex eigenvector, respectively. The subscripts $s$ and $k$ denote the symmetric part and skew symmetric part. For any matrix $M$, the associated $M_s = \frac{M+M^T}{2}$ and $M_k=\frac{M-M^T}{2}$. Write $x_i=x_{Ri}+jx_{Ii}$ with $x_{Ri}$ and $x_{Ii}$ being real and imaginary parts, respectively. Multiplying the above equation by $x_i^*$ on both sides yields the following second-order scalar equation,
\begin{equation} \label{III12}
(a_{mi}+jb_{mi})\lambda_i^2+(a_{di}+jb_{di})\lambda_i+(a_{ki}+jb_{ki})
 =0, \quad i=1,\cdots,2n,
\end{equation} 
where $x_i^*$ is the conjugated transpose of $x_i$, $a_{mi}=x_{Ri}^TM_sx_{Ri}+x_{Ii}^TM_sx_{Ii}$, $b_{mi}=2x_{Ri}^TM_kx_{Ii}$, and $a_{di}$, $b_{di}$, $a_{ki}$, $b_{ki}$ are similarly expressed. The solutions of \eqref{III12} satisfy the following result,
\begin{lemma} \label{th04}
The solution $\lambda_i$ of \eqref{III12} has negative real parts if and only if 
\begin{equation} \label{III13}
b_{mi}b_{di}+a_{mi}a_{di}>0,
\end{equation}
and 
\begin{equation} \label{III14}
(a_{di}a_{ki}+b_{di}b_{ki})(a_{mi}a_{di}+b_{mi}b_{di})>(a_{mi}b_{ki}-b_{mi}a_{ki})^2.
\end{equation}
\end{lemma} 
\bigskip

The above lemma is closely related to the symmetric degree of matrices. If all the involved matrices are symmetric positive definite, then \eqref{III13} and \eqref{III14} hold and the corresponding dynamic system is asymptotically stable. 

To describe in which degree a matrix $M$ is symmetric positive definite, we define a  measurement variable as follows
\begin{equation} \label{eqa02}
\theta(M)=\max_{\theta} \left\{\theta\ge 0\left| \begin{bmatrix} M_s & 0 \\ 0 & M_s \end{bmatrix}\pm \theta \begin{bmatrix}
 0 & M_k \\ M_k^T & 0 \end{bmatrix}\ge 0 \right.\right\}
\end{equation}

The above matrix inequality implies that $M_s\ge 0$ and so does $M$. If $M$ is negative definite, then the matrix inequality in \eqref{eqa02} is infeasible and therefore $\theta(M)$ does not exist. 

The physical meaning of the above definition is that complex matrix $M_s+j\theta M_k$ formed by the symmetric parts and asymmetric parts of $M$ is positive semi-definite for all $0\le\theta\le\theta(M)$. If $M$ is symmetric, then $\theta(M)\rightarrow \infty$, which means that no asymmetric parts exist; if the symmetric parts of $M$ is zero, then $\theta(M)=0$, which means that $M$ is  {skew symmetric}.

\begin{theorem} \label{th05}
System \eqref{III04} is semistable if the network is connected and one of the following conditions is satisfied,
\begin{description}
\item [c1).] there are positive scalars $\bar{u}$ and $\bar{r}$ such that $R_iI_i^{max}=\bar{u}$ and $\alpha_i=\beta_i\bar{r}$ for all $i\in\mathcal{N}$, a positive scalar $\nu_1$ such that 
\begin{equation} \label{III16}
\Upsilon+\bar{u}Y\ge \nu_1 Y\Psi L
\end{equation}
and $\theta(Y\Psi L)$ exists such that
\begin{equation} \label{III17}
(\nu_1+\bar{r})\theta^2(Y\Psi L)+\bar{r}>\tau_{max},
\end{equation}
where $\tau_{max}=\max \{\tau_1,\cdots,\tau_n\}$ is the maximum time constant of nodal low pass filter. 
\item [c2).] there are positive scalars $\tau$, $\bar{r}$, $\beta$ and $\nu_2$ such that $D=\tau E$, $\Phi = \bar{r}\Psi$, $\Psi = \beta E$,
\begin{equation} \label{III18}
Y^{-1}\Upsilon+R\Upsilon\ge \nu_2 \beta L,
\end{equation}
$Y_s\ge 0$, and $\theta(Y\Upsilon)$ exists such that
\begin{equation} \label{III19}
(\nu_2+\bar{r})\theta^2(Y\Upsilon)+\nu_2+\bar{r}>\tau.
\end{equation} 
\end{description}
\end{theorem}
\begin{IEEEproof}
We firstly prove that system \eqref{III04} has only one zero eigenvalue when the network is connected. Rewrite $A_c$ by 
\begin{equation}
A_c=A_{c1}A_{c2}=A_{c3}A_{c4},
\end{equation}
with
\begin{gather*}
A_{c1}=\begin{bmatrix} D^{-1} & 0 \\ -\Phi L (\Upsilon)^{-1}D^{-1} & E \end{bmatrix}, \quad
A_{c2}=\begin{bmatrix} -(E+YR) & Y \\ -\Psi L (\Upsilon)^{-1} & 0 \end{bmatrix}, \\
A_{c3}=\begin{bmatrix} E & 0 \\ -\Phi L (\Upsilon)^{-1} & - \Psi L (\Upsilon)^{-1}\end{bmatrix}, \\
A_{c4}=\begin{bmatrix} -D^{-1}(E+YR) & D^{-1}Y \\ E & 0\end{bmatrix}.
\end{gather*}
The rank of system matrix $A_c$ equals to that of $A_{c2}$, as well as to that of $A_{c3}$. For a connected network, $rank(L)=rank(Y_c)=n-1$ and both them have $\underline{1}_n^T$ and $\underline{1}$ as the associated left and right eigenvector respectively to the trivial eigenvalue $0$. It is known that if $Y_s\neq 0$, namely at least one shunt conductance exists, then $Y=Y_c+Y_s>0$. Below we proceed by two cases.

i) $Y_s\neq 0$. In this case $Y$ is a M-matrix whose inverse matrix $Y^{-1}$ is nonnegative matrix. It can be seen that $\rank(A_{c2})=n+\rank(L)=2n-1$ due to $Y$ being nonsingular. Therefore $0$ is one eigenvalue of $A_c$. We further show that the $0$ is semisimple by showing $\mathrm{rank}(A_c^2)=2n-1$ (according to Proposition 1 in the appendix). Noting that $\mathrm{rank}(A_{c}^2)=\mathrm{rank}(A_{c2}A_{c3})$, we consider the null space of $A_{c2}A_{c3}$. Suppose there is a vector $x=[x_1^T,x_2^T]^T\in\real^{2n}$ such that
\begin{equation}\label{III21}
-A_{c2}A_{c3}x=
\begin{bmatrix} E+YR+Y\Phi L(\Upsilon)^{-1} & Y\Psi L(\Upsilon)^{-1} \\
\Psi L (\Upsilon)^{-1} & 0 \end{bmatrix}
\begin{bmatrix} x_1 \\ x_2 \end{bmatrix}=0
\end{equation}
The solution of $x$ has $x_1=\sigma_1 \Upsilon\underline{1}$ for some $\sigma_1\in\real$, which as well should satisfy
\begin{equation}
\sigma (Y^{-1}+R)\Upsilon\underline{1}=\Psi L (\Upsilon)^{-1}x_2.
\end{equation}
Left multiplying $\underline{1}_n^T\Psi^{-1}$ on both sides leads to 
\begin{equation}
\sigma_1 [\beta_1,\beta_2,\cdots,\beta_n](Y^{-1}+R) [I_1^{max}, I_2^{max}, \cdots, \Upsilon_n] =0,
\end{equation}
which implies $\sigma_1=0$ because positive definite symmetric matrix $Y^{-1}+R$ is nonnegative. Thus, the solution of \eqref{III19} has the form of $\sigma_2[0,\Upsilon\underline{1}_n^T]^T$ for some $\sigma_2\in\real$. This means that $\mathrm{rank}(A_{c2}A_{c3})=2n-1$, so dose $A_c^2$. Subsequently $A_c$ has only one zero eigenvalue. 

ii) $Y_s=0$. In this case $Y=Y_c$ is singular. It can be verified that all the solutions of $A_{c2}x=0$ has the form of $x=\sigma_1[0, \Upsilon \underline{1}_n^T]^T$ for some $\sigma_1\in\real$, which means that $\mathrm{rank}(A_c)=\mathrm{rank}(A_{c2})=2n-1$. Below we further show $\mathrm{rank}(A_c^2)=2n-1$. Suppose there is a vector $x=[x_1^T,x_2^T]^T\in\real^{2n}$ such that $A_c^2x=0$, then $A_cx=\sigma_1[0, \Upsilon \underline{1}_n^T]^T$ for some scalar $\sigma_1$, which implies 
\begin{equation}
A_{c2}x=\begin{bmatrix} -(E+Y_cR) & Y_c \\ -\Psi L (\Upsilon)^{-1} & 0 \end{bmatrix}
\begin{bmatrix} x_1 \\ x_2 \end{bmatrix}
=\sigma_1 A_{c1}^{-1}\begin{bmatrix} 0 \\ \Upsilon\underline{1}_n \end{bmatrix}
=\begin{bmatrix} 0 \\ \sigma \Upsilon\underline{1}_n \end{bmatrix}.
\end{equation}
The second block line gives $-\Psi L (\Upsilon)^{-1}x_1=\sigma_1\Upsilon\underline{1}_n$.  Left multiplying $\underline{1}_n^T\Psi^{-1}$ on both sides obtains $0=\sigma_1\underline{1}_n^T \Psi^{-1}\Upsilon\underline{1}_n$, which implies $\sigma_1=0$. Thus, it can be concluded that if $A_c^2{x}=0$ then $A_cx=0$. That means $\mathrm{rank}(A_c^2)=\rank(A_c)=2n-1$ and subsequently $A_c$ has one simple zero eigenvalue as well.

The reminder is to prove that all nonzero eigenvalues are of negative real parts for every condition.

c1) With the condition c1), recasting system \eqref{III04} into the form of \eqref{III11} yields the following matrix parameters,
\begin{gather}
M_s=D\Upsilon, D_s=\Upsilon+\bar{u}Y+\bar{r}K_s, D_k=\bar{r}K_k\\
 K_s=\frac{Y\Psi L+L\Psi Y}{2},\quad K_k = \frac{Y\Psi L - L\Psi Y}{2}.
\end{gather}
Due to the existence of $\theta(Y\Psi L)$, $K_s\ge 0$ and subsequently $D_s>0$. Therefore $a_{mi}>0$ and $a_{di}>0$ for all $i$. If $a_{ki}=0$, then $b_{ki}=0$ due to $a^2_{ki}\ge \theta^2(Y\Psi L)b^2_{ki}$, and the nonzero solution of \eqref{III12} is $-\frac{a_{di}+jb_{di}}{a_{mi}}$, being of negative real parts.

Now consider the case of $a_{ki}>0$. Because of $b_{mi=0}$, \eqref{III13} holds always and \eqref{III14} reduces to
\begin{equation}
(a_{di}a_{ki}+b_{di}b_{ki})a_{di}>a_{mi}b_{ki}^2.
\end{equation}
If $b_{ki}=0$ and $a_{ki}\neq 0$, then the above is obvious. 
Matrix inequality \eqref{III16} means $a_{di}\ge (\nu_1+\bar{r})a_{ki}$. On the other hand  $a_{di}>\frac{1}{\tau_{max}}a_{mi}$ because of $a_{ki}>0$.  Noticing $b_{di}=\bar{r}b_{ki}$ and $a^2_{ki}\ge \theta^2(Y\Psi L)b_{ki}$, one has 
\begin{equation} \label{III28}
(a_{di}a_{ki}+b_{di}b_{ki})a_{di}
>(\nu_1+\bar{r})\theta^2(Y\Psi L)+\bar{r})b^2_{ki}\frac{a_{mi}}{\tau_{max}},
\end{equation}
which together with \eqref{III17} and $b_{ki}\neq 0$ implies \eqref{III28} and therefore \eqref{III14}. By lemma \ref{th04}, all eigenvalues associated with $a_{ki}> 0$ have negative real parts. 

Now It can be concluded that system \eqref{III14} has only one simple $0$ eigenvalues and all other eigenvalues be of negative real parts. Thus it is semistable.

\medskip
c2) With condition c2), $Y>0$. Multiplying by $Y^{-1}$ on both sides of \eqref{III05} and then casting it into the form of \eqref{III11} obtain
\begin{gather}
M_s=\tau \frac{Y^{-1}\Upsilon+\Upsilon Y^{-1}}{2}, \quad D_s=R\Upsilon+\bar{r}K_s+M_s/\tau, \\ K_s=\beta L,\quad
M_k = \tau \frac{Y^{-1}\Upsilon-\Upsilon Y^{-1}}{2}, \quad D_k=M_k/\tau.
\end{gather} 
Since $\theta(Y^{-1}\Upsilon)=\theta(Y\Upsilon)$ exists, $a_{mi}\ge 0$ and $a_{mi}^2\ge \theta^2(Y\Upsilon) b_{mi}^2$. Also $a_{di}>\frac{a_{mi}}{\tau}$ and $a_{di}>(\nu_2+\bar{r})a_{ki}$ due to \eqref{III18}. 

Firstly consider $b_{mi}=0$. If $a_{mi}=0$ as well, then \eqref{III12} has the solution $\lambda_i=\frac{-a_{ki}(a_{di}-jb_{di})}{a^2_{di}+b^2_{di}}$ that is either $0$  or of negative real parts when $a_{ki}=0$ or $a_{ki}>0$, respectively. If $a_{mi}>0$ but $a_{ki}=0$, then similarly it can be seen that the nonzero solution of \eqref{III12} has negative real parts. 

Now consider $b_{mi}\neq 0$. Noting that $b_{di}=\frac{b_{mi}}{\tau}$, \eqref{III13} holds always, which as well ensures that the nonzero solution of \eqref{III12} with $a_{ki}=0$ has negative real parts. For $a_{ki}\neq 0$, \eqref{III14} reduces to 
\begin{equation} \label{III31}
a_{di}a_{ki}(a_{mi}a_{di}+b_{di}b_{mi})>b^2_{mi}a^2_{k_i}.
\end{equation}
Noticing that $a_{di}a_{ki}\ge (\nu_2+\bar{r})a_{ki}^2$ and $a_{mi}a_{di}+b_{di}b_{mi}>(\theta^2(Y\Upsilon)+1)\frac{b^2_{mi}}{\tau}$, inequality \eqref{III13} holds by \eqref{III19}. This by lemma \ref{th04} shows that all the eigenvalues associated with $b_{mi}\neq 0$ and $a_{ki}>0$ have negative real parts. Therefore system \eqref{III14} has only one simple $0$ eigenvalue and all other eigenvalues be of negative real parts and is semistable.
\end{IEEEproof}
\medskip

Notice that $\bar{r}$ can be set a value larger than $\tau_{max}$ that in general is a small value less than $0.1$. Therefore conditions from \eqref{III16} to \eqref{III19} are easy to satisfy. But more critical are the implied conditions in theorem \ref{th05}, the existence of $\theta(Y\Psi L)$ and $\theta(Y\Upsilon)$ for cases c1) and c2), respectively. If $I^{max}_1=I^{max}_2=\cdots=I_n^{max}$, then $\theta(Y\Upsilon)$ always exists. This means the heterogeneous extent between nodes influences the stability of DC microgrids with distributed control \eqref{III02}.

It should be pointed out that the information network and physical network are not required to have the same topology. For a special case that $L=Y_c$ and moreover $Y_s=gE$ and $\Psi=\beta E$, then $\theta(Y\Psi L)$ always exists as well.   

\subsection{Steady state}
The steady state of closed-loop system \eqref{III04} is the mode determined by the $0$ eigenvalue. This subsection issues not only what is the steady state but also the relationship with the initial condition and network topology. 

\begin{theorem}
A semistable dynamic system \eqref{III04} will reach the current sharing in the sense that all nodes have the same current ratio $i_k^r=r_{c1}$ for all $k\in\mathcal{N}$, where $r_{c1}$ is given by,
\begin{equation} \label{III32}
r_{c1}=\left\{\begin{array}{ll}
 \dfrac{\underline{1}_n^T\Psi^{-1}\big(U^d(0)+\Phi L (\Upsilon)^{-1}I^{m}(0)\big)}{\underline{1}_n^T\Psi^{-1}(Y^{-1}+R)\Upsilon\underline{1}_n},  & \mbox{if } Y_s\neq 0 \\
0 & \mbox{if } Y_s=0
\end{array}
\right..
\end{equation}
Moreover, the output voltage of nodes is given by
\begin{equation} \label{III33}
U=\left\{\begin{array}{ll}
r_{c1}Y^{-1}\Upsilon\underline{1}_n & \mbox{if } Y_s \neq 0 \\
r_{c2}\underline{1}_n & \mbox{if } Y_s = 0
\end{array}
\right.,
\end{equation}
where $r_{c2}$ is defined by
\begin{equation} \label{III34}
r_{c2}
=\frac{\underline{1}_n^T\Psi^{-1}(U^d(0)+\Phi L (\Upsilon)^{-1}I^{m}(0))}{\underline{1}_n^T\Psi^{-1}\underline{1}_n}.
\end{equation}
\end{theorem}
\begin{IEEEproof}
With $Y_s\neq 0$,the left eigenvector $v_l$ and the right eigenvector $v_r$ of $A_c$ associated with $0$ eigenvalue are respectively,
\begin{equation}
v_l=\begin{bmatrix} (\Upsilon)^{-1}L\Phi \Psi^{-1}\underline{1}_n  \\ \Psi^{-1}\underline{1}_n \end{bmatrix}, \quad \quad v_r = \begin{bmatrix} \Upsilon\underline{1}_n \\ (Y^{-1}+R)\Upsilon\underline{1}_n \end{bmatrix}. 
\end{equation}
The left eigenvector $v_l$ leads to 
\begin{equation}
\frac{d}{dt}v^T_l\begin{bmatrix} I^m(t) \\ U^d(t) \end{bmatrix}\equiv 0;
\end{equation}
while the semistable implies that the trajectory of \eqref{III04} converges to the space spanned by the right eigenvector $v_r$, that is, there is a scalar $r_{c2}$ such that 
\begin{equation}
\lim_{t\rightarrow \infty}
\begin{bmatrix} I^m(t) \\ U^d(t) \end{bmatrix}=r_{c2} v_r.
\end{equation}
Noting $L\underline{1}_n=0$, a combination of the above two formulae yields 
\begin{equation}
r_{c2}=\frac{v^T_l[(I^m(0))^T, (U^d(0))^T]^T}{v_l^Tv_r} 
=\dfrac{\underline{1}_n^T\Psi^{-1}\big(U^d(0)+\Phi L (\Upsilon)^{-1}I^{m}(0)\big)}{\underline{1}_n^T\Psi^{-1}(Y^{-1}+R)\Upsilon\underline{1}_n}.
\end{equation}

Below consider the case of $Y_s=0$ for which $Y=Y_c$ is singular. The left and right eigenvector of $A_c$ associated with $0$ eigenvalue are respectively
\begin{equation}
v_l = \begin{bmatrix} (\Upsilon)^{-1}L\Phi \Psi^{-1}\underline{1}_n  \\ \Psi^{-1}\underline{1}_n \end{bmatrix}
\quad \mathrm{and}\quad 
v_r = \begin{bmatrix} 0 \\ \underline{1}_n \end{bmatrix}.
\end{equation}
Similarly, it can be given that 
\begin{equation}
r_{c2}=\frac{v^T_l[(I^m(0))^T, (U^d(0))^T]^T}{v_l^Tv_r}
=\dfrac{\underline{1}_n^T\Psi^{-1}\big(U^d(0)+\Phi L (\Upsilon)^{-1}I^{m}(0)\big)}{\underline{1}_n^T\Psi^{-1}\underline{1}_n},
\end{equation}
and subsequently $\lim_{t\rightarrow \infty} I^m(t)=0$ and $\lim_{t\rightarrow\infty} U^d(t)=r_{c2}\underline{1}_n$. 

Equation \eqref{III33} follows from the fact that $U(t)=U^d(t)-RI^{m}(t)$. This completes the proof. 
\end{IEEEproof}

In general the initial values $I^{m}(0)$ of the low-pass filter for output currents are set to zeros; or else they will influence the steady states according to \eqref{III32}, which is not what we want. In the following, we always assume that $I^{m}(0)=0$. With this, the proportional gain $\alpha_i$ will influence the stability but not the steady state.

It is rational to assume that the steady currents in the decentralized droop control, stated in \eqref{eqII12}, do not excess the maximum permissible currents for all nodes. That is, 
\begin{equation} \label{III41}
(E+YR)^{-1}YU^d(0)\preceq I^{max}, 
\end{equation}
where $I^{max}\in\real^{n}=[I^{max}_{1},\cdots, I_n^{max}]^T$ and $\prec$ denotes the  { component} less than. Define the current ratio vector by $R_{cr}=[i_1^r,\cdots,i_n^r]^T$. Denote the steady current ratio in the decentralized droop control by $R^{Dec}_{cr}$, which satisfies $R^{Dec}_{cr}=\Upsilon^{-1}(E+YR)^{-1}YU^{d}(0)$. Denote by $\overline{R}^{Dec}_{rc}$ and $\underline{R}^{Dec}_{rc}$ the maximum and minimum elements of $R^{Dec}_{cr}$, respectively.

With $I^{m}(0)=0$, several further discussions on the steady states for some special cases are made. 
\begin{corollary} \label{th07}
Given $I^{m}(0)=0$ and $Y_s\neq 0$, the steady states of system \eqref{III04} with the cooperative droop control satisfy the following properties:
\begin{enumerate}
\item [P1)] With \eqref{III41}, the common current ratio satisfies $r_{c1}<1$ and  $r_{c1}\in [\underline{R}_{cr}^{Dec},\overline{R}_{cr}^{Dec}]$.
\item [P2)] With \eqref{III41}, the output voltage satisfies $\frac{r_{c1}}{\overline{R}_{cr}^{Dec}}U^{ss}\preceq U\preceq \frac{r_{c1}}{\underline{R}_{cr}^{Dec}}U^{ss}$.
\item [P3)] If $\Psi=\beta E$ and $R=rE$, then 
$r_{c1}<\frac{\sum_k U_k^d(0)}{\sum_k I_k^{max}}\dfrac{\lambda_n}{1+\lambda_n r}$, where $\lambda_n$ is the maximal eigenvalue of $Y$.
\item [P4)]  
If $Y_s=gE$, $\Psi=\beta E$, and $R=rE$, then $r_{c1}=\frac{\sum_k U_k^d(0)}{\sum_k I_k^{max}}\dfrac{g}{1+gr}$.
\item [P5)] All the nodes have the same output voltage if and only if there is a positive scalar $\epsilon$ such that $\Upsilon=\epsilon Y_s$.
\end{enumerate}
\end{corollary}
\begin{IEEEproof} P1)
When $Y_s\neq 0$, $Y$ is a M-matrix whose inverse matrix  $Y^{-1}$ is a nonnegative matrix. Note that $\Upsilon \underline{1}_n=I^{max}$ and \eqref{III41}, it follows that
\begin{equation}
r_{c1}<\frac{\underline{1}_n^T\Psi^{-1}U^d(0)}{\underline{1}_n^T\Psi^{-1}(Y^{-1}+R)(E+YR)^{-1}YU^d(0)}=1.
\end{equation}
Also noting that $U^{d}(0)=(Y^{-1}+R)\Upsilon R^{Dec}_{cr}$, one has 
\begin{equation} \label{III43}
r_{c1}=\frac{\underline{1}_n^T\Psi^{-1}(Y^{-1}+R)\Upsilon R^{Dec}_{cr}}{\underline{1}_n^T\Psi^{-1}(Y^{-1}+R)\Upsilon \underline{1}_n}.
\end{equation}
Since $\underline{1}_n^T\Psi^{-1}(Y^{-1}+R)\Upsilon\succ 0$, 
\begin{equation}
\underline{1}_n^T\Psi^{-1}(Y^{-1}+R)\Upsilon \underline{R}^{Dec}_{cr}\underline{1}_n^T\preceq \Psi^{-1}(Y^{-1}+R)\Upsilon R^{Dec}_{cr}
 \preceq \underline{1}_n^T\Psi^{-1}(Y^{-1}+R)\Upsilon \bar{R}^{Dec}_{cr}\underline{1}_n,
\end{equation}
which together with \eqref{III43} shows $r_{c1}\in[\underline{R}^{Dec}_{cr},\bar{R}^{Dec}_{cr}]$.

\smallskip
P2) Recall the steady output voltage under decentralized droop control,
\begin{equation}
U^{ss}=(E+RY)^{-1}U^d(0)=(E+RY)^{-1}(Y^{-1}+R)\Upsilon R^{Dec}_{cr}
 =Y^{-1}\Upsilon R^{Dec}_{cr}.
\end{equation}
Due to the nonnegativeness of $Y^{-1}$, $\underline{R}_{cr}^{Dec}Y^{-1}I^{Max}\preceq U^{ss}\preceq \bar{R}_{cr}^{Dec}Y^{-1}I^{Max}$. This together with $U=r_{c1}Y^{-1}I^{Max}$ shows the property P2).

\smallskip
P3) Noticing that $Y^{-1}>\frac{1}{\lambda_n}E$ and $Y^{-1}\succ 0$, $\underline{1}_n^TY^{-1}\succ \frac{1}{\lambda_n}1_n^T$. With this, one has $r_{c1}<\frac{\sum_k U_k^d(0)}{\sum_k I_k^{max}}\dfrac{\lambda_n}{1+\lambda_nr}$.

\smallskip
P4) Recall the proof of lemma \ref{th02}, there is a transformation matrix $T$ and a diagonal matrix $\Gamma = \diag(0,\lambda_2,\cdots,\lambda_n)$ such that $Y_c = T^{-1}\Gamma T$ where $T^{-1}=[\underline{1}_n,\star]$ and $T=[\underline{1}_n\frac{1}{n}, \star]^T$. With this, $Y^{-1}+R=T^{-1}\big((gE+\Gamma)^{-1}+rE\big)T$. Noting that $\underline{1}_n^TT^{-1}=[n,0,\cdots,0]$, it can be derived that $r_{c1}=\frac{\sum_k U_k^d(0)}{\sum_k I_k^{max}}\dfrac{g}{1+gr}$.

\smallskip
P5) It suffices to show that $Y^{-1}I^{max}=\epsilon \underline{1}_n$ for some scalar $\epsilon$. This, together with $Y_c\underline{1}_n=0$, means that $\Upsilon=\epsilon Y_s$.
\end{IEEEproof}   

\medskip
Three remarks are presented in order. 
\begin{remark}
P1) means that the same current ratio under the cooperative droop control locates between the minimal current ratio and the maximal current ratio under the decentralized droop control. The node with the minimal current ratio will increase its output current; on the contrary the output current of the node with the maximal current ratio will decrease.  
\end{remark}

\begin{remark}
 {In the case that all the nodes have the same controller gains and the same virtual resistance}, P3) shows an upper bound of current ratio related to the maximal eigenvalue of conductance matrix $Y$. Furthermore, if the  {shunt} conductances of nodes are the same  as well, then P4) shows that $r_{c1}$ equals to the ratio between the current of shunt conductance $g$ under the average voltage and the average maximal current.  {Both of them are novel but limit to the special cases}.   
\end{remark}

\begin{remark}
P5) shows that in general the sharing of both current and voltage can not occur simultaneously.  {This point has also been stated in \cite{LuTPE2014}.}
\end{remark}

\begin{figure}[htbp]
\begin{center}
\includegraphics[width=0.4\textwidth]{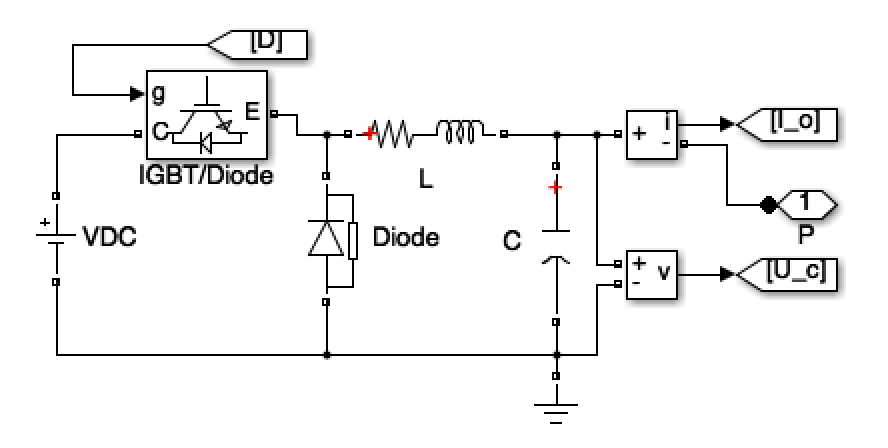}
\includegraphics[width=0.5\textwidth]{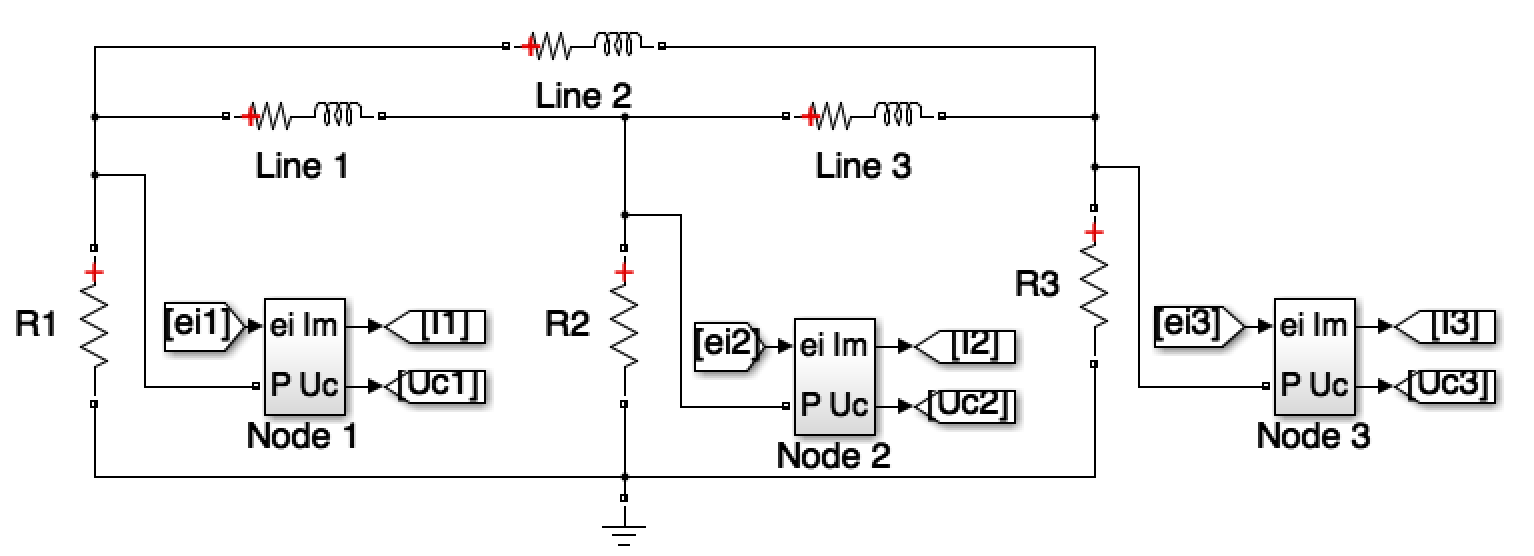}
\caption{Three nodes DC microgrid. The buck converter (top) and network structure (bottom).}
\label{fig01}
\end{center}
\end{figure}

\section{Simulation example} \label{sec04}
Consider a microgrid with three nodes connected by a triangular form. Each node is a DC voltage interfaced to the DC bus via a buck converter. As illustrated in Fig. \ref{fig01}, the simulation is made on the Matlab Simulink. The electrical parameters are
listed in Table \ref{tab01}.
\begin{table}[htp]
\caption{Electrical Parameters}
\begin{center}
\begin{tabular}{c|c|c|c|c|c}
\hline
Line 1 & Line 2 & Line 3 & R1 & R2 & R3 \\
\hline
$1\Omega$, $0.5$mH & $0.5\Omega$, $0.1$mH & $0.4\Omega$, $0.2$mH & $2\Omega$ & $5\Omega$
& $4\Omega$\\
\hline 
\end{tabular}
\end{center}
\label{tab01}
\end{table}%

The same control parameters are selected for all the nodes. The desired DC bus voltage is $48$V, the maximum current of all nodes is $30$A, the time constant of LPF is $0.01$s, the virtual resistance $r=0.1\Omega$, the proportional gain $\alpha=0$ and the integrator gain $\beta=100$. The conductance matrix of network and the Laplacian matrix of information graph are given by, respectively,  
\begin{equation*}
Y = \begin{bmatrix}3.5 & -1 & -2 \\ -1 & 3.7 & -2.5 \\ -2 & -2.5 & 4.75 \end{bmatrix}, \quad 
L = \begin{bmatrix}2 & -1 & -1 \\ -1 & 2 & -1 \\ -1 & -1 & 2 \end{bmatrix}.
\end{equation*}
It can be verified that $YL$ is not a positive definite matrix, therefore $\theta(Y\Psi L)$ does not exist. Since $\Upsilon$ is an unit matrix, the conditions \eqref{III18} and \eqref{III19} hold and subsequently the microgrid with distributed controller \eqref{III02} is semistable. 

\begin{figure}[htbp]
\begin{center}
\includegraphics[width=0.49\textwidth]{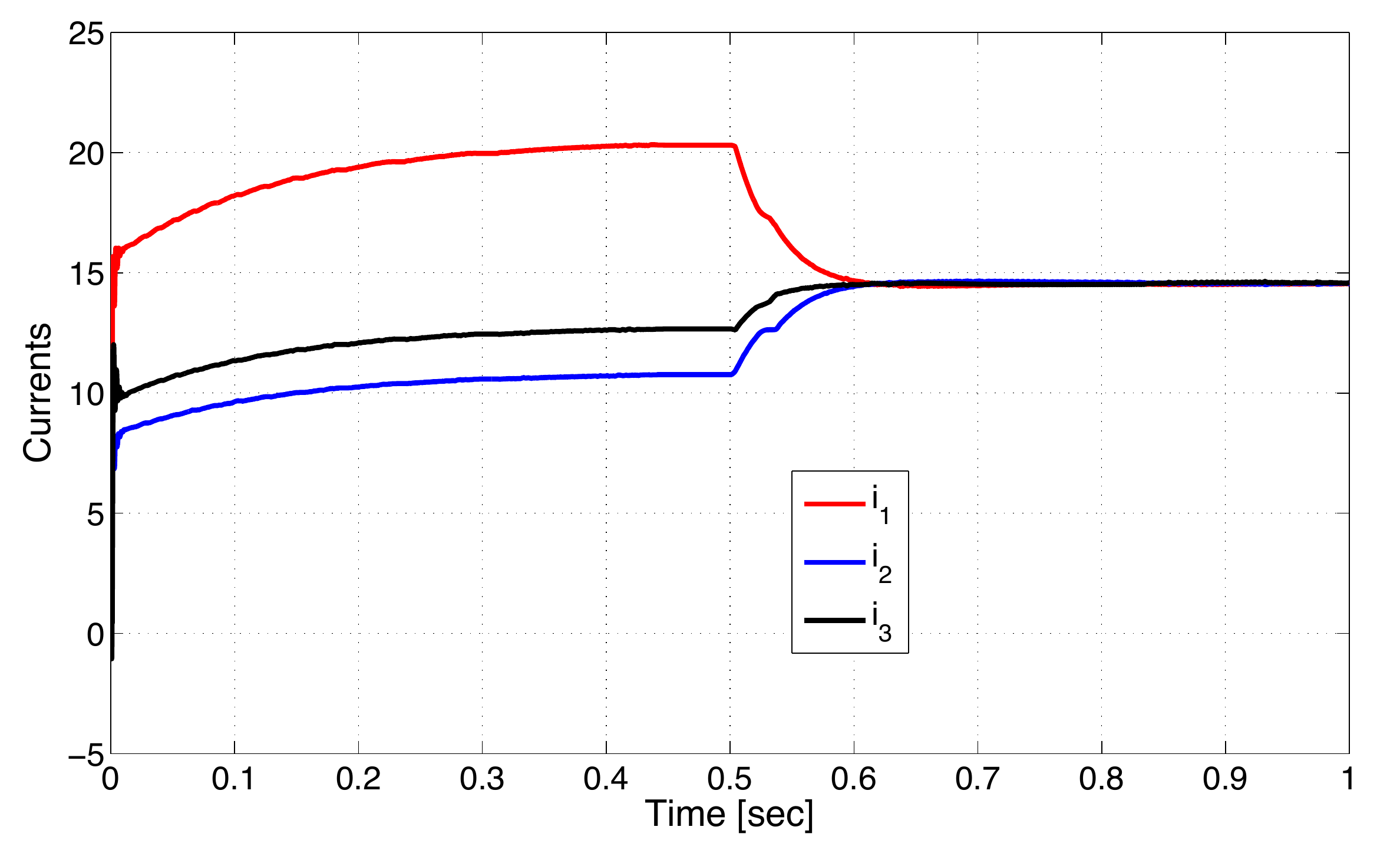}
\includegraphics[width=0.49\textwidth]{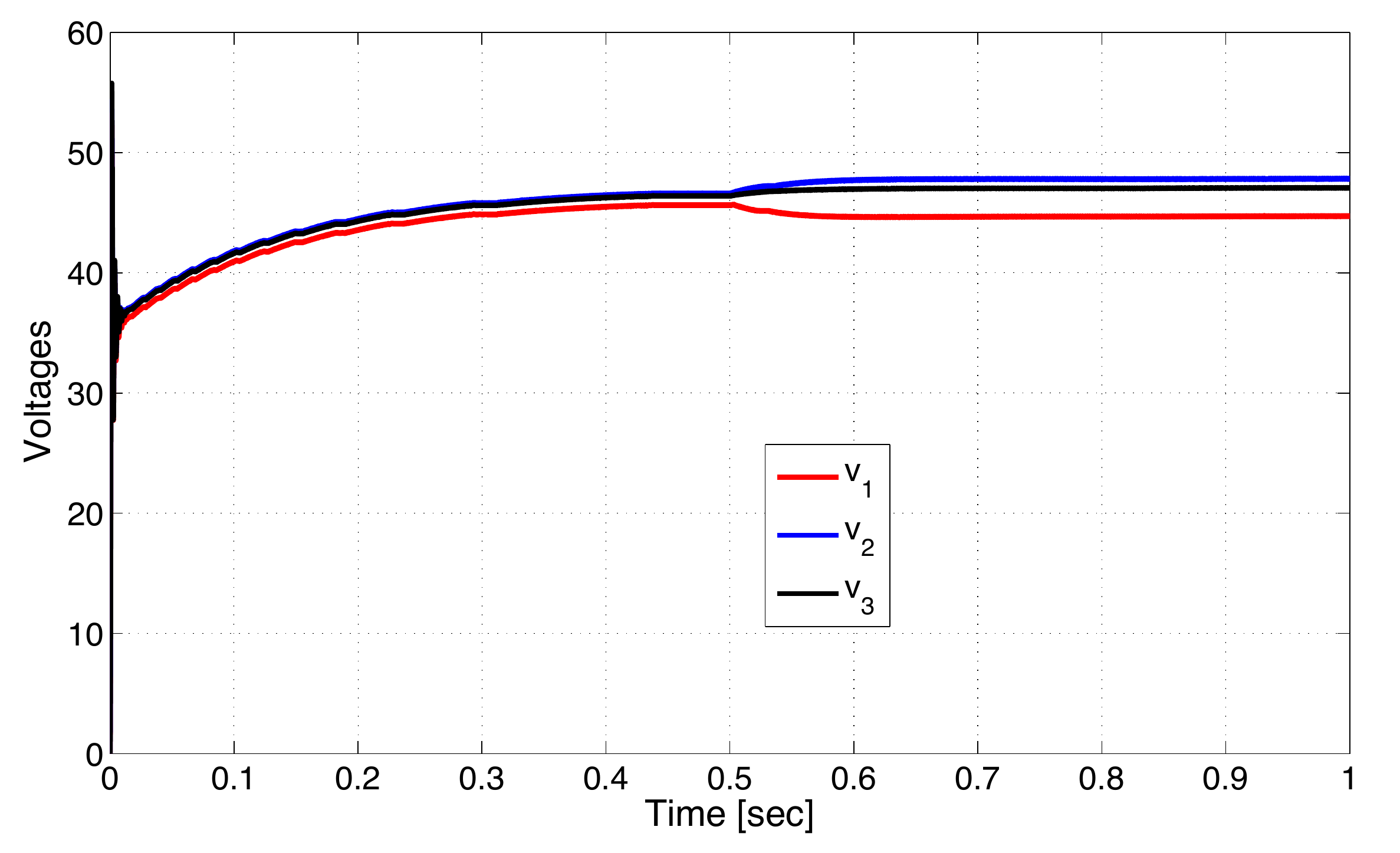}
\caption{Trajectories of currents and voltages of three nodes}
\label{fig02}
\end{center}
\end{figure}

Fig. \ref{fig02} shows the simulation results. At the beginning, the microgrid runs in the decentralized droop control. At $t=0.5$s, the system is closed to the steady state under the decentralized droop control, which are
\begin{equation}
I^{ss}=\begin{bmatrix}20.31\\10.76\\12.66 \end{bmatrix},\quad U^{ss}=\begin{bmatrix} 45.64 \\ 46.60 \\ 46.41 \end{bmatrix}
\end{equation}
In the presences of line resistance and the difference of local load, there is a large bias for current sharing. The minimal and maximal current ratio are $\underline{R}_{cr}^{Dec}=0.36$ and $\overline{R}_{cr}^{Dec}=0.68$.

At the time $t = 0.5$s, the distributed cooperative control is applied, with which the nodal currents then asymptotically converge to the same value, an exact current sharing. All the nodes have the same output current $14.56$A and current ratio $0.4853$ that belongs to $[\underline{R}_{cr}^{Dec},\overline{R}_{cr}^{Dec}]$. The steady output voltages are $U=[44.71,47.83,47.07]^T$ satisfying property P2 in corollary \ref{th07}.

\section{Conclusion} \label{sec05}
The stability of DC microgrids with distributed cooperative control has been investigated and two sufficient semistable conditions are presented. The study on steady state illustrated the current sharing property and its relationship with initial condition and network topology. A DC microgrid with three buck converter nodes was simulated on the Matlab platform to show the developed results.

The developed results based on the constant impedance matrix $Y$ are applicable as well for the time-varying $Y$ with a duly small fluctuation. The stability with a general time-varying $Y$ is one of our ongoing research, for which the analysis through the locations of eigenvalues as used in this paper is not applicable any more.

\appendix
Given a linear system $\dot{x}(t)=Ax(t)$ where $t>0$, $x(t)\in\real^n$ and $A\in\real^{n\times n}$, the following definitions are made \cite{BernsteinJMD1995},\\
\textbf{Definition 1.} The system is {\em Lyapunov stable} if, for every initial condition $x(0)$, there exists $\epsilon>0$ such that $\|x(t)\|<\epsilon$ for all $t\ge 0$.\\
\textbf{Definition 2.} The system is {\em semistable} if $\lim_{t\rightarrow \infty} x(t)$ exists for all initial conditions $x(0)$.\\
\textbf{Definition 3.} Given an eigenvalue $\lambda\in \mbox{spec}(A)$, $\lambda$ is {\em semisimple} if every Jordan block of $A$ associated with $\lambda$ is of size one. Further, it can be seen that $\lambda$ is semisimple if and only if 
\begin{equation}
\mathrm{rank}(\lambda I-A)=rank(\lambda I-A)^2,
\end{equation}
where $I$ is the unit matrix.\\

The following proposition is true. \\
\textbf{Proposition 1.} $A$ is semistable if and only if $A$ is Lyapunov stable and $A$ has no nonzero imaginary eigenvalues.

\bibliographystyle{unsrt}
\bibliography{/Users/yanjun/JabRef/XJ2013}

\end{document}